\documentclass[12pt]{article}
\usepackage{amsfonts,amsthm}
 \usepackage[cp1251]{inputenc}

\usepackage[T2A]{fontenc}
\usepackage[english]{babel}
\theoremstyle{plain} \tolerance 9000 \hbadness 9000

\def\R{{\mathbb R}}
\def\C{{\mathbb C}}
\begin{document}

\title{
    \bf
Entire functions of exponential type, almost periodic in Besicovitch's sense on the
real hyperplane}

\author{S.Yu.~Favorov, O.I.~Udododva}

\date{}

\maketitle

{\it 2000 Mathematics Subject Classification:} {\small Primary 32A15, Secondary 42A75,
32A50}

{\it Keywords:} {\small  entire function of exponential type, Besicovitch's almost
periodic function}

\newtheorem{Th}{Theorem}
\newtheorem*{Le}{Lemma}
\newtheorem*{ThA}{Theorem A}
\newtheorem*{ThB}{Theorem B}

\begin{abstract}
{ Suppose that an almost periodic in Besicovitch's sense function $f(x)$ of several
variables is the restriction to the real hyperplane of an entire function of
exponential type $b$. Then its spectrum is contained in the ball of radius $b$ with
the center in the origin. }
\end{abstract}

In his paper \cite{Bohr} H.~Bohr showed that the spectrum of an almost periodic
function $f(x)$ on the real axis $\R$ is a subset of $[-\sigma, \sigma]$, as long as
$f$ is the restriction to $\R$ of an entire function of an exponential type $\sigma$.
 R.~Boas \cite{Boas} extended the assertion to
almost periodic functions  on  $\R$ in Besicovitch's metric (for brevity, B-almost
periodic functions). In the general case, these functions are unbounded on $\R$, hence
the proof of the latter assertion is more difficult.

H.~Bohr's result was generalized  to almost periodic functions in a finite dimensional
space by S.Yu.~Favorov and O.I.~Udodova \cite{FavUd}. But the case of B--almost
periodicity is more complicated, because restrictions to straight lines of B--almost
periodic functions  in  $\R^p$ are not necessary almost periodic.

It should be mentioned that B--almost periodic functions of several variables
 were considered earlier in  \cite{Gir}, \cite{Ud}, \cite{U}.
But in \cite{Ud}, \cite{U} the spectrum of functions was not under consideration, and
in \cite{Gir} the author studied only B--almost periodic functions with bounded
Besicovitch's norm in a tube domain with a cone in the base.

Our proof differs from ones in \cite{Boas}, \cite{Bohr}, \cite{FavUd} and is based on
estimates of entire functions and Logvinenko's theorem \cite{Logvinenko} on the growth
of entire functions of several variables on the hyperplane $\R^p$.

We will use the following notations.

By $z=x+iy,$ $z=(z_1,...,z_p),$ $x=(x_1,...,x_p),$ $y=(y_1,...,y_p)$ we denote the
vectors in $\C^p$ (or, respectively, in $\R^p$, $'x$ means the vector $(x_2,...,x_p)
\in \R^{p-1}$, $\langle x, y \rangle$ is the inner product in $\R^p$. Next, $|z|$,
$|x|$, $|'x|$ are the Euclidean norms in the spaces $\C^p$, $\R^p$, and $\R^{p-1}$,
respectively. By   $dx$, $d'x$, and  $dx_1$ we denote the Lebesgue measure in $\R^p$,
$\R^{p-1}$, and $\R$, respectively. Furthermore, $B(x, \delta)$ means the open ball in
$\R^p$ of radius $\delta$ with the center in $x$, $C$ with lower indexes are
constants, depending only on $f$.

{\it Besicovitch's norm} of a locally integrable function $f(x)$ in $\R^p$ is the
limit
 $$
  \|f\|_B=
\mathop {\overline{\lim}}\limits_{T \to \infty } \left( {\frac{1}{{2T}}} \right)^p
\int\limits_{[-T,T]^p}^{}  \,|f(x)|dx.
 $$
The function $f(x)$ is called {\it B--almost periodic} in $\R^p$, if for any
$\varepsilon
> 0$ there is a (generalized) trigonometric polynomial
\begin{equation}\label{polinom}
 P(x)=\sum c_n e^{i\langle x, \lambda^{(n)}
\rangle}, \quad c_n \in {\bf C}, \,
 \lambda^{(n)} \in {\bf R^p},
\end{equation}
 such that
 $$
 \|f-P\|_B < \varepsilon.
 $$
{\it The Fourier coefficient} of $f$ is the limit
 \begin{equation} \label{a1}
 a(\lambda , f)= \mathop {\lim}\limits_{T \to \infty } \left(
{\frac{1}{{2T}}} \right)^p \int\limits_{[-T,T]^p}^{} \,f(x)e^{-i \langle x, \lambda
\rangle}dx.
 \end{equation}
{\it The spectrum} $\texttt{sp} f$ of $f(x)$ is the set
$$
 \{
\lambda \in {\bf R^p} : a(\lambda, f)\neq 0\}.
 $$
 Note that the existence the limit (\ref{a1}) and countability of the spectrum follow easily
from the definition of B--almost periodicity, and the  equality
\begin{equation}
\label{a2}
a\left( {\lambda ,P} \right) = \left\{ \begin{array}{l}
 c_n ,\;\lambda  = \lambda ^{\left( n \right)}  \\
 0,\;\lambda  \ne \lambda ^{\left( n \right)},  \\
 \end{array}\right.
\end{equation}
which holds for any $n$ and any polynomial (\ref{polinom}).

By \cite{U},  for any B--almost periodic function there exists a sequence of
polynomial (\ref{polinom}) (the so-called Bochner-Fejer sums) with
$\lambda^{n}\in\texttt{sp} f$, which approximate $f$.

The main result of our paper is the following theorem.

\begin{Th}\label{1}
 Let a B--almost periodic function $f(x)$ in $\R^p$ extends to $\C^p$
as an entire function with the bound
\begin{equation}
\label{f}
 |f(z)| \leq C_0 e^{\sigma |z|}.
\end{equation}
 Then we have
  ${\rm{sp}} f \subset
B(0, \sigma)$.
\end{Th}

The proof of this theorem is based on the following statement.

\begin{Th}\label{2}
Let $f(x), x \in \R^p$ be a function with a finite norm  $\|f\|_B$. If $f$ can be
extended to $\C^p$ as an entire function  with estimate (\ref{f}), then
 \begin{equation}
\label{prod} |f(x)| \leq C_1 \prod\limits_{j = 1}^p {(1 + |x_j |)^p } \qquad\forall x
\in \R^p.
\end{equation}
 \end{Th}
We get the proof of theorem \ref{2}, using the following auxiliary results.
\begin{ThA}[\cite{Logvinenko}]
Let  $f(z)$ be an entire function on $\C^p$, which satisfies (\ref{f}), and $E$ be a
$\delta$-net in $\R^p$. If $\sigma \delta < K(p)$, then
 $$
\mathop {\sup }\limits_{x \in R^p } |f(x)| \le \left( {1 - \sigma \delta } \right)^{ -
1} \mathop {\sup }\limits_{x \in E} |f(x)|.
 $$
\end{ThA}
 \begin{ThB}[see, for example, \cite{Levin}, p. 311]
Let a function $g(w)$ be a holomorphic  in $\C^+=\{w \in \C, {\rm Im}\,w>0\}$,
continuous in the closure of $\C^+$, and satisfy the estimate
 \begin{equation}
\label{g(w)} |g(w)| \leq ce^{a|w|}, \, w \in \C^+.
\end{equation}
If
\begin{equation}
\label{int} \int\limits_{ - \infty }^{ + \infty } {\frac{{\log ^ + |g(t)|}}{{1 + t^2
}}dt}  < \infty,
\end{equation}
then
$$\log |g(w)| \le \frac{{{\mathop{\rm Im}\nolimits} w}}{\pi
}\int\limits_{ - \infty }^{ + \infty } {\frac{{\log |g({\mathop{\rm Re}\nolimits} w +
t)|}}{{t^2  + ({\mathop{\rm Im}\nolimits} w)^2 }}dt}  + h{\mathop{\rm Im}\nolimits}
w,\;\;w \in \C^+,
 $$
 where $h=\mathop {\overline {\lim } }\limits_{t \to + \infty }
\frac{{\log |g(it)|}}{t}.$
 \end{ThB}
In the case  $\sup_\R |g(w)|<\infty$ theorem B yields the well known version of
Fragment--Lindel\"{o}f Principle
\begin{equation}
\label{gFL}
 |g(w)| \leq
 \mathop {\sup }\limits_{{\mathop{\rm Im}\nolimits}\, w = 0}
 |g(w)|e^{h{\mathop{\rm Im}\nolimits}\, w} ,\;\;w \in\C^+.
\end{equation}
{\bf Proof of Theorem \ref{2}.}
 Since $\|f\|_B < \infty$, we get for any $\delta \in
(0,1)$ and $\widetilde x \in \R^p$
$$\int\limits_{B(\widetilde x,\delta )} {|f(x)|dx}
\le \int\limits_{[-\widetilde x-\delta , \widetilde x+\delta ]^p}
 {|f(x)|dx}\leq C_2\left( {1 + |\widetilde x|} \right)^p. $$
Therefore there is a constant $C_3<\infty$ such that for any  ball of radius $\delta$
there exists a point $x'$ in the ball with
$$|f(x')| \le {{C_3 }}{{\delta ^{-p} }}\left( {1 + |x'|} \right)^p. $$
Put
 $$
 g(z) = f(z)\prod\limits_{j = 1}^p {\left( {\frac{{\sin z_j }}{{z_j }}} \right)^p } .
 $$
We have
 $$
 |g(z)| \le C_0 e^{\left( {\sigma  + p^2 } \right)|z|},\quad z \in \C^p.
 $$
Since $|g(x)| \le C_4\delta^{-p}$ at the points of the  $\delta$-net, we see that
theorem A with a suitable  $\delta$  implies the bound
\begin{equation}
\label{sup} \mathop {\sup }\limits_{x \in \R^p} |g(x)| \le C_5.
\end{equation}
Using (\ref{sup}), we apply inequality (\ref{gFL}) first in the domain ${\rm Im} z_1
> 0$, and then in the domain ${\rm Im} z_1 < 0$. We get
\begin{equation}
\label{g'}
 |g(z_1,\, 'x)| \le C_5 e^{\left( {\sigma + p^2 }
\right)|y_1 |}
\end{equation}
 for all $z_1=x_1+iy_1 \in \C$,  $'x \in \R^{p-1}$. Apply
 (\ref{gFL}) to the function
 $g(z_1,\dots,z_p)$ as a function in the variable $z_2$ and use (\ref{g'}) instead of (\ref{sup}).
Repeating these arguments by the variables $z_3, \ldots , z_p$, we get
 $$
 |g(z)| \le c_5 e^{\left( {\sigma  + p^2 } \right)\left( {|y_1 |
+  \ldots  + |y_p |} \right)}.
 $$
 Therefore,
 $$
 \left| {f(z)\cdot \prod\limits_{j = 1}^p {\left( {\frac{{\sin z_j }}{{z_j }}}
\right)^p } } \right| \le C_6
 $$
 on the set $A=\{ z \in {\bf C^p} : \, |y_1| \le 1,
\ldots |y_p| \le 1\}$. Hence,
 $$\left| {f(z)\cdot \prod\limits_{j = 2}^p {\left( {\frac{{\sin z_j }}{{z_j }}}
\right)^p } } \right| \le C_7(1+|z_1|)^p
 $$
 on the set $\{z: z\in A, z\notin
\bigcup \limits_{n} B(n\pi, \frac{1}{2})\}$ . By the Maximum Principle, we get the
same inequality with the constant $2^p C_7$ instead of $C_7$ at every point of $A$.
Repeating these arguments  $p-1$ times, we obtain (\ref{prod}). Theorem is proved.

\bigskip
For the proof of Theorem \ref{1} we need the following Lemma.
\begin{Le}
Suppose that $f(z)$ is an entire function in $\C^p$, which satisfies (\ref{f}), and
its restriction to $\R^p$ satisfies the condition $\|f\|_B \le \infty$. Then for any
$s_0 \in (0, \infty)$,  $T \geq T(s_0)$, and $s \in (0, s_0)$ we get
$$\int\limits_{[ - T,T]^p } |f(x_1+is,'x)|e^{-s \sigma}dx \le C_8 T^p,$$
with $C_8=2^{p+1}(1+2\cdot 3^p\|f\|_B)$.
\end{Le}
{\bf Proof of the Lemma.}
 By theorem \ref{2}, the function $f(z_1,\, 'x)$
satisfies (\ref{int}) in the variable $z_1$ for any fixed  $'x \in \R^{p-1}$. Taking
into account (\ref{f}), we get
 $$
 \overline {\mathop {\lim }\limits_{t \to \infty } } \frac{{\log
|f(iy_1 ,\;'x)|}}{{y_1 }} \le \sigma.
 $$
 Hence Theorem B implies for any $s> 0$
 $$
\log |f(x_1  + is,\;'x)| \le \frac{s}{\pi }\int\limits_{ - \infty }^{ + \infty }
{\frac{{\log |f(t + x_1 ,\;'x)|dt}}{{t^2  + s^2 }}} + \sigma s.
 $$
Since the measure $\frac{s}{\pi}\cdot\frac{dt}{t^2+s^2}$ is a probability one on
 $\R$, we get for any locally integrable function $h(t)$ on $\R$
 $$
 \exp\left({\frac{s}{\pi}
\int\limits_{|t|\le 2T}\frac{h(t)dt}{t^2+s^2}}\right)\le
\frac{s}{\pi}\int\limits_{|t|\le
2T}\frac{e^{h(t)}dt}{t^2+s^2}+\frac{s}{\pi}\int\limits_{|t| > 2T}\frac{dt}{t^2+s^2}.
 $$

Next,
$$\int\limits_{[ - T,T]^p } {|f(x_1  +
is,\;'x)|e^{ - \sigma s}dx }  \le \int\limits_{[ - T,T]^p } {\exp \left\{
{\frac{s}{\pi }\int\limits_{ - \infty }^{ + \infty } {\frac{{\log |f(t + x_1
,\;'x)|}}{{t^2  + s^2 }}dt} } \right\}dx}
$$
$$\le\int\limits_{[ - T,T]^p } {\exp \left\{ {\frac{s}{\pi }
\int\limits_{|t| \le 2T}^{} {\frac{{\log |f(t + x_1 ,\;'x)|}}{{t^2  + s^2 }}dt} }
 \right\}
\times }$$ $$\exp \left\{ {\frac{s}{\pi }\int\limits_{|t| > 2T}^{}
  {\frac{{\log |f(t + x_1 ,\;'x)|}}{{t^2  + s^2 }}dt} } \right\}dx \leq $$
$$\int\limits_{[ - T,T]^p } {\left( { \frac{s}{\pi }
\int\limits_{|t| \le 2T}^{} {\frac{{|f(t + x_1 ,\;'x)|}}{{t^2  + s^2 }}dt}
 }+1 \right)\exp \left\{ {\frac{s}{\pi }\int\limits_{|t| > 2T}^{}
 {\frac{{\log |f(t + x_1 ,\;'x)|}}{{t^2  + s^2 }}dt} } \right\}dx}.$$

By (\ref{prod})  we have for  $x \in \left[ { - T,\;T} \right]^p $
$$\frac{s}{\pi }\int\limits_{|t| \ge 2T}^{}
{\frac{{\log |f(t + x_1 ,\;'x)|}}{{t^2  + s^2 }}dt} \le \frac{s}{\pi }\int\limits_{|t|
\ge 2T} {\frac{{C_3\left(\sum\limits_{j=1}^{p}p\log^+ (1+|x_j|)+p\log^+t\right)}}{{t^2
+ s^2 }}dt}
$$
 $$
 \le C_3\left[ p^2\log(1+T) \frac{s}{\pi}\int\limits_{|t|\ge 2T}\frac{dt}{t^2
+s^2}+\frac {ps}{\pi}\int\limits_{|t|\ge 2T}\frac{\log^+tdt}{t^2+s^2}\right].
 $$
Note that the latter expression bounds from above by $\log 2$ for $s \le s_0$ and $T
\ge T(s_0)$.

Therefore, taking into account the inequality
 $$
 \int\limits_{[ - T,T]^p } {|f(x)|dx} \le 2\|f\|_B \left( {2T} \right)^p,\quad T \ge C_{9}
 $$
we obtain
$$\int\limits_{[ - T,T]^p } {|f(x_1  + is,\;'x)|e^{ - s\sigma } dx}
\le 2\left( { \frac{s}{\pi }\int\limits_{|t| \le 2T} {\frac{{dt}}{{t^2 + s^2
}}\int\limits_{[ - T,T]^p } {|f(x_1  + t,\;'x)|dx} } }+(2T)^p \right)$$
$$
\le 2^{p + 1} T^p  + 2\frac{s}{\pi }\int\limits_{|t| \le 2T} {\frac{{dt}}{{t^2  + s^2
}}} \left[ {\int\limits_{[ - 3T,3T]} {|f(x)|dx} } \right] \le 2^{p+1}(1+2\cdot 3^p
\|f\|_B) T^p. $$
 Lemma is proved.

\medskip
 {\bf Proof of Theorem \ref{2}}

Let $A$ be an orthogonal matrix in $\R^p$. It is easy to check the equality
\begin{equation}
\label{alambda11}
 a(\lambda, f)=a(A^{-1} \lambda, f_A),
\end{equation}
where $f_A(x)=f(Ax)$ and $\lambda$ is an arbitrary vector in $\R^p$.

Indeed, it follows from (\ref{a2}) that this equality is true for any polynomial $P$
(\ref{polinom}). To prove it for an arbitrary B--almost periodic function, we can
approximate it by polynomial $P$ such that $\|f-P\|_B~<~\varepsilon$. Therefore,
$\|f_A-P_A\|_B < K^p \varepsilon$, where $K=\mathop {\max }\limits_j |Ae_j |$, $e_j$
is the natural basis in
 $\R^p$. Hence, we obtain (\ref{alambda11}).

Take $\lambda \in \R^p$, $|\lambda| > \sigma$. Since bound  (\ref{f}) is the same for
$f_A$, we may suppose that $\lambda = (-\sigma-\eta, 0, \ldots , 0),$ $\eta > 0.$ In
this case we have
\begin{equation}
\label{alambda12}
 a(\lambda, f)=\mathop {\lim }\limits_{T \to \infty } \left( {\frac{1}{{2T}}}
  \right)^p \int\limits_{[ - T,T]^p } {f(x_1 , \, 'x)e^{ix_1 (\sigma  + \eta )}
  dx}.
\end{equation}
The function $f(z_1, x')$ is holomorphic in $z_1$, therefore for any $y_1>0$ we have

$$\int\limits_{[ - T,T]^p } {f(x_1, \, 'x)e^{ix_1 (\sigma  + \eta )}
dx}=  $$
$$ i\int\limits_0^{y_1 } {\int\limits_{[ - T,T]^{p - 1} }
{f( - T + is,\, 'x)e^{ - iT(\sigma  + \eta ) - s(\sigma  + \eta )}
 d\,'x~ds} } +$$
$$
\int\limits_{[ - T,T]^p } {f(x_1  + iy_1 , \, 'x)e^{ix_1 (\sigma +
\eta ) - y_1 (\sigma  + \eta )} dx} - $$
$$ i\int\limits_0^{y_1
} {\int\limits_{[ - T,T]^{p - 1} } {f(T + is,\;'x)e^{iT(\sigma  + \eta ) - s(\sigma  +
\eta )} d\;'x~ds} } $$ $$=I_1 (T, y_1)+ I_2(T, y_1)-I_3(T, y_1).
$$

By Lemma, we get
$$|I_2(T,y_1)|\le C_8 T^pe^{-\eta y_1},$$
hence for a given $\varepsilon >0$ and sufficiently large  $y_1$
\begin{equation}
\label{lim} \overline{\lim\limits_{T\rightarrow\infty}}|(2T)^{-p}I_2(T,y_1)|\le
\varepsilon.
\end{equation}
Next, $|I_1-I_3| \le G(T)$, where
$$
G(x_1)=\int\limits_0^{y_1 } {\int\limits_{[ - T,T]^{p - 1} } {e^{-s\sigma}\left(|f(x_1
+ is,\;'x)| + |f( - x_1 + is,\;'x)|\right)~ds~d\;'x}}.
$$
By Lemma, the Lebesgue measure of the set $$E=\left\{ {x_1 :\frac{T}{2} < |x_1 | <
T,\;G(x_1 )
>3C_8 T^{p - 1} |y_1 |} \right\}$$ is at most
$$\frac{1}{{2C_8 |y_1 |T^{p - 1} }}\int\limits_E {G(x_1 )dx_1 }\leq$$
$$
 \frac{1}{{3C_8 |y_1 |T^{p - 1} }}\int\limits_0^{y_1 }
{\int\limits_{[ - T,T]^p } { e^{-s\sigma}\left( {|f(x_1  +
is,\;'x)| +
 |f( -
x_1  + is,\;'x)|} \right)dx~ds}}<\frac{2T}{3}.
 $$
  Hence for some  $T' \in \left[ {\frac{T}{2},T}
\right]\backslash E$ we get $G(T')\leq 3 C_8 |y_1| T^{p-1}$. Therefore we have
$$\overline{\lim\limits_{T'\rightarrow \infty}}\left( {2T'} \right)^{ - p}
\left| { (I_1(T',y_1)-I_3(T',y_1)) } \right| =0.
 $$
Thus, the latter bound and (\ref{lim}) yield $a(\lambda, f)=0$. Theorem is proved.

\bigskip
Mathematical School, Kharkov national university, Svobody sq. 4,
Kharkov, 61077, Ukraine.

e-mail: Sergey.Ju.Favorov@univer.kharkov.ua

\bigskip

 Department of Mathematics, Ukrainian State Academy of Railway Transport,
Feyerbah sq. 7, Kharkov, 61050, Ukraine.

e-mail: udodova@kart.edu.ua

\end{document}